\documentclass[numbook, envcountsect, envcountsame, envcountreset, runningheads, smallextended]{svjour3}

\usepackage{amscd,amsfonts,amssymb,amsmath,latexsym,array,hhline,xcolor,graphicx}

\usepackage{latexsym}
\usepackage{times}

\newcommand\cE{{\cal E}}

\newcommand\cF{{\cal F}}

\newcommand\cL{{\cal L}}

\newcommand\e{{\varepsilon}}

\def\E{{\bf E}}
\def\P{{\bf P}}

\def\bbr{{\mathbb R}}

\def\text#1{\hbox{#1}}

\def\E{{\bf E}}
\def\P{{\bf P}}

\def\build #1_#2{\mathrel{\mathop{\kern 0pt #1}\limits_{#2}}} 
\newcommand{\zs}[1]{{\mathchoice{#1}{#1}{\lower.25ex\hbox{$\scriptstyle#1$}}
{\lower0.25ex\hbox{$\scriptscriptstyle#1$}}}}

\numberwithin{equation}{section}

\def\bbr{{\mathbb R}}

\def\bbr{{\mathbb R}}
\newcommand\fdem{$\Box$}
\newcommand\beq{\begin{equation}}
\newcommand\eeq{\end{equation}}
\newcommand\bea{\begin{eqnarray}}
\newcommand\eea{\end{eqnarray}}
\newcommand\bean{\begin{eqnarray*}}
\newcommand\eean{\end{eqnarray*}}

\newtheorem{theo}{Theorem}[section]

\newtheorem{assm}[theo]{Assumption}

\begin{document}

	\title{Distributional equations and the ruin problem for the Sparre Andersen model with investments 
	}
	\author{ Yuri Kabanov \and Danil Legenkiy \and Platon Promyslov} 
	
	\institute{\at	Lomonosov Moscow State University, Moscow, 119991, Russia
	and Universit\'e de Franche-Comt\'e, Laboratoire de Math\'ematiques, UMR CNRS 6623, 
	16 Route de Gray, 25030 Besan\c{c}on, France \\
 \email{ykabanov@univ-fcomte.fr},    
    \and \at Lomonosov Moscow State University and ``Vega'' Institute Foundation, Moscow, 119991, Russia
    \\ 
    \email{danil.legenky@gmail.com},
	\and \at Lomonosov Moscow State University and ``Vega'' Institute Foundation, Moscow, 119991, Russia\\
	 \email{platon.promyslov@gmail.com}}

\titlerunning{Distributional equations and the ruin probabilities for a Sparre Andersen model with investments}

\date{\today}
\maketitle

\begin{abstract}
This note is an addendum to the work initiated by Eberlein, Kabanov, and Schmidt and developed further by Kabanov and Promyslov on the asymptotics of the ruin probabilities in the Sparre Andersen model with investments in a risky asset. Using more advanced methods of the implicit renewal theory, we provide complements to some results of the mentioned works.
\end{abstract}

\keywords{
Ruin probabilities \and Sparre Andersen model \and Actuarial models with investments \and Renewal processes \and Annuities \and Distributional equations 
}
\smallskip

\noindent
 {\bf Mathematics Subject Classification (2010)} 60G44
 
 \medskip
\noindent
 {\bf JEL Classification} G22 $\cdot$ G23 
 
\section{Introduction}
In modern world insurance companies usually place their reserves, fully or partially, in financial markets. This practice motivated a development of mathematical theory of ruin with risky investments whose origin can be traced back to the pioneering work \cite{Paul-93} by Paulsen who developed it further in \cite{Paul,Paul-98,Paul-G}; The main tool was the Kesten--Goldie theory of distributional equations as exposed in \cite{Go91}. The recent progress in the latter led to essential improvements of asymptotical results on ruin, see \cite{KP2020,KPro,KPukh}. Our note is inspired by a comment of an anonymous referee of \cite{KPro} who suggested to use deeper results from the theory of distributional equations to simplify proofs and obtain weaker sufficient conditions of the mentioned paper. Following this advice, we use here a corrolary following from Proposition 2.5.4 from the fundamental book by Buraczewski, Damek, and Mikosch, \cite{BDM}, and get new sufficient conditions for the property of interest. As we already mentioned, models with risky investments depend on characteristics of the involved processes. Though the asymptotic analysis is based on a few general principles already exploited in previous papers, the new tool and related estimates, specific for every considered variant, require a certain mathematical skill and may be of interest for the reader working in the ruin theory with investments or in the theory of distributional equations and its applications. 
 
Recall that there is an alternative approach to the ruin problems with investments based on the asymptotic analysis of integro-differential equations of the second order for the ruin probabilities, see  \cite{ACT,BelkinaKK,Fr,FrKP} {\it et al.}
In this short note we leave aside a detailed discussion of the evolution of this theory and send the reader to the recent papers \cite{AK}, \cite{KP2020}. 
We mention here only the main conclusion: under reasonable assumptions, the ruin probability in models with risky investments decays as a power function of the initial capital in a striking contrast with the exponential decay in the classical Cram\'er--Lundberg theory. Of course, this fact motivates further research in this field.

The ruin theory with investments is a vast domain: models of the insurance business are combined with numerous models for the price processes of risky assets or stochastic interest rates. Study of various combinations is far from being a routine exercise. 
Here we consider the combination where the business is described by the Sparre Andersen model, that is, by a compound renewal process with drift, while the price of the risk asset follows an independent L\'evy process. This setting is not a new one: e.g., the asymptotic of ruin probabilities in a version of the annuity payments (upward jumps of the business process) was studied by Eberlein, Kabanov, and Schmidt, \cite{EKS}, the non-life insurance version (downward jumps) and a mixture of both were treated by Kabanov and Promyslov, \cite{KPro}. Both papers assume that the law of claims has light tails and the interarrival times have an exponential moment. The principal results assert that the ruin probabilities, for any version of the business process, admit a power function asymptotic if the price process is of unbounded variation or its L\'evy measure charges both half-axes. The looking simpler cases where the price process is of bounded variation with only positive or only negative jumps happen to be 
rather delicate. The main difficulty is to check whether the support of the solution of an affine distributional equation is unbounded from above. It can be formulated as a property of a limit at infinity of specific stochastic processes, see, e.g., a paper by Bohme {\it et al.} \cite{BLRR} by Bohme {\it et al.} 

The results established in \cite{KPro} involve additional assumptions on the supports of laws of claims and interarrival times. The natural question arises whether one can do better. 
The approach mentioned above allows us to consider the cases not covered in the previous studies.

The paper uses standard notations of stochastic calculus and concepts discussed in details in \cite{EKS,KP2020}.

\section{The model}
We consider the setting where the process 
$X=(X_t)_{t \ge 0}$ is given as the solution of non-homogeneous linear stochastic equation 
\begin{equation}
\label{X}
dX_t = X_{t-} dR_t + dP_t, \quad X_0=u>0, 
\end{equation}
where $P$ is a compound renewal process independent of the L\'evy process $R$ with $\Delta R>-1$. It covers several important models of the collective risk theory where 
$X$ describes the capital reserve evolution. 

The case where $R$ is equal to zero, i.e. $X_t=u+P_t$, is just the classical Sparre Andersen model (see, e.g. \cite{Gr}, Ch.3) where the capital reserve $X$ evolves only due to the business activity process $P$. The latter is usually represented by the formula
\begin{equation*}
P_t=ct+ \sum_{i=1}^{N_t}\xi_i,
\end{equation*}
where $N=(N_t)_{t \ge 0}$ is a renewal counting process with jumps at moments $T_i$, $i \ge 0$. ``Renewal'' means that the interarrival times $T_i-T_{i-1}$ form an i.i.d. sequence independent on the i.i.d. sequence $(\xi_i)$. By assumption, 
$T_0: = 0$. The common distribution function of the random variables $\xi_i=\Delta P_{T_i}$ is denoted by $F_\xi$ (the same notation is often used also for the law $
\cL(\xi)$). The common distribution function of the interarrival times $T_i-T_{i-1}$, $i\ge 2$, is denoted by $F$. 

\smallskip
The theory distinguishes the two major variants of the model: 

\smallskip
$(a)$ Non-life insurance: $c>0$, $\xi_i<0$. The company receives the incoming stream of payments from its customers and covers claims of the insured. 

\smallskip
$(b)$ Annuity payments: $c<0$, $\xi_i>0$. The company pays pensions and gets back the remains 
of customers deposits. In a more recent literature this case is interpreted as a model of the business activity of a venture company paying salary and selling innovations. Of course, for the latter interpretation the use of exponentially distributed interarrival times as in the Cram\'er--Lundberg model seems not to be adequate.

Less frequently, in the literature, see, e.g. \cite{ACT}, one can find the third variant, interesting mainly from a mathematical point of view and corresponding to the case where the company combines two types of activity:

\smallskip
$(c)$ The random variables $\xi_i$ take values in $\bbr\setminus \{0\}$. 

\smallskip
To understand the general model given by (\ref{X}) suppose that the company invest instantaneously its total capital in a risky asset with the prices $S=(S_t)$ which is a geometric L\'evy process with an infinitesimal relative increment $dR_t=S_{t-}/dS_t$. In other words, $S$ is the solution of the linear stochastic equation $dS_t=S_{t-}dR_t$, that is, $S=S_0\cE(R)$ where $\cE(R)$ is the stochastic exponential of $R$. 
We can rewrite the first term in the right-hand side of (\ref{X}) as the product of $X_{t-}/S_{t-}$ (the number of units of the risky asset) and $dS_t$ (the price increment). That is, the first term describes the capital reserve increment due to the risky investment while $dP_t$ describe the increment due to the business activity.

Let $(a,\sigma ^ 2,\Pi)$ be the L\'evy triplet of $R$. We assume that $R$ is not deterministic, that is at least one of the conditions $\sigma^2 \neq 0$ or $\Pi \neq 0$ is fulfilled.
The condition $\Delta R>-1$ or, in the equivalent 
form $\Pi ((-\infty,-1])=0$, ensures that 
the price process $S>0$. The stochastic exponential $S=S_0\cE(R)$ can be written as the usual exponential $e^V$. The logprice process $V$ happens to be the L\'evy process with the triplet 
$(a_V,\sigma^ 2,\Pi_V)$, where
$$
a_V=a-\frac{\sigma^2}{2}+\Pi \big (h( \ln(1+x))-h\big ), 
$$
with $h: x\mapsto xI_{\{|x|\le 1\}}$, $x\in {\mathbb R }$, and $\Pi_V=\Pi\varphi^{-1}$, $\varphi: x\mapsto \ln (1+x)$, $x> -1$. Here and in the sequel 
a function in the argument of a measure means its integral with respect to this measure. In the standard notations of stochastic calculus the process $V$ can be represented by the formula
\begin{equation*} V_t=at-\frac 12 \sigma^2 t + \sigma W_t+ h*(\mu-\nu)_t+\big (\ln (1+x)-h\big )*\mu_t,
\end{equation*}

We define the moment of ruin as $\tau^u:=\inf \{t\colon X^u_t\le 0\}$ and consider the ruin probability as a function of the initial capital by putting $\Psi(u):={\bf P}[\tau^u<\infty]$.

We exclude from consideration the cases $c\ge 0$, $\xi>0$ (the ruin never happens) and $c<0$, $\xi\le 0$ (the ruin happens for sure). 

Let $H:q\mapsto \ln {\bf E} [e^{-qV_{T_1}}]$ be the cumulant generating function of the random variable $V_{T_1}$. 

We shall work under two standing assumptions. 

\begin{assm} 
\label{ass1}
The function $H$ has a root $\beta>0$ such that $H(\beta+)\neq +\infty$.
\end{assm}

\begin{assm} 
\label{ass2}
${\bf E}[|\xi_1|^{\beta}]<\infty$ and ${\bf E}[e^{\varepsilon T_1}]<\infty$ for some $\varepsilon>0$. 
\end{assm}

The aim of our note is to use the recent development in the implicit renewal theory to relax some sufficient conditions of papers \cite{EKS} and \cite{KPro} ensuring the properties 
\begin{equation}
\label{maina}
0<\liminf_{u\to \infty}u^{\beta}\Psi(u)\le \limsup_{u\to \infty}u^{\beta}\Psi(u)<\infty. 
\end{equation}

In the sequel we shall use the following notations: 
$$
\underline{F}_\xi := \inf \operatorname{supp} F_\xi,\ \ \overline{F}_\xi:=\sup \operatorname{supp} F_\xi,
\quad \underline F:= \inf \operatorname{supp} F, \ \ \overline{F}:=\sup \operatorname{supp} F,\qquad
$$
and, in the same spirit, $ \underline{\Pi}_V := \inf \operatorname{supp} \Pi_V, \quad \overline{\Pi}_V:=\sup \operatorname{supp} \Pi_V$, etc. 

Note that $[\underline{F}_\xi, \overline{F}_\xi]$ is the convex hull of the support of law $F_\xi$. 
\smallskip

We reformulate the results of the mentioned papers in fairly transparent form as two theorems. The first theorem covers all three versions 
 of models: non-life insurance, annuity payments and the mixed one, cf. Ths 1.1 - 1.3 in \cite{KPro}. 

\begin{theorem} 
\label{old1}
Suppose that $\sigma^2>0$, or $\Pi( |h|)=\infty$, or 
$\Pi$ charges $(-1,0)$ and $(0,\infty)$, or $\underline F_\xi=-\infty$. Then (\ref{maina}) holds. 
\end{theorem}

In particular, if $R$ (or $S$) is of unbounded variation then the first two conditions are satisfied (see Prop. 3.9 in \cite{CT}), hence, (\ref{maina}) holds. 

\smallskip

\begin{theorem} 
\label{old2} The property (\ref{maina}) holds if $\Pi$ charges either $(-1,0)$ or $(0,\infty)$ and, 

-- for the non-life insurance model $\underline F=0$, 

-- for the annuity payment model $\overline{F}=\infty$,

-- for the mixed model $\overline{F}_{|\xi|}=\infty$ or $\overline{F}=\infty$, if $c<0$, and $\underline F_{\xi^+}=0$, if $c\ge 0$. 
\end{theorem}

We could formulate the above theorem with a disjoint hypothesis by adding the negation of the assumption of Theorem \ref{old1}, that is, $\sigma^2=0$, $0<\Pi(|h|)<\infty$, and $\underline F_\xi>-\infty$. 
In the following theorems, which are the main results of our note, 
we also will tacitly assume this negation which implies, in particular, that the logprice process has the form $V=(a-\Pi (h))t+x*\mu^V$. In the case, where $\Pi( (-\infty, 0) ) = 0$, that is the logprice has only positive jumps, the difference $a - \Pi(h)<0$ (otherwise Assumption \ref{ass1} will be violated). Since it is convenient to work with positive parameters we put $\lambda:=\Pi (h)-a$. In the case where the logprice has only negative jumps, i.e. $\Pi(0,\infty)=0$, we put $\kappa:=a-\Pi (h)$ and this value is strictly positive to avoid a contradiction with Assumption \ref{ass1}.

\begin{theorem}
\label{new1}
Suppose that $\Pi$ charges either $(-1,0)$ or $(0,\infty)$. Then (\ref{maina}) holds:

-- for the non-life insurance model when $\underline F>0$, if $\Pi((-\infty,0))=0$ and \\ $ \overline{F}_{|\xi|}>(c/\lambda) (1-e^{-\lambda \underline F})$
or $\Pi((0,\infty)) = 0$ and $\overline{F}_{|\xi|}>(c/\kappa) (e^{\kappa \underline F}-1) $, 

-- for the annuity payment model when $\overline{F}<\infty$, if $\Pi((-\infty,0))=0$ and \\ $\underline F_{\xi}<(c/\lambda) (1-e^{-\lambda \overline{F}})$
or $\Pi((0,\infty)) = 0$ and $\underline F_{\xi}<(c/\kappa)(e^{\kappa \overline{F}}-1) $.
\end{theorem}

\begin{theorem}[mixed case]
\label{new2}
Suppose that $\Pi$ charges either $(-1,0)$ or $(0,\infty)$.  Then (\ref{maina}) holds:

-- when $c\! \ge\! 0$ and $\underline{F}>0$, if 
$\Pi((-1, 0)) = 0$ and $\overline{F}_{\xi^-}>(c/\lambda) (1-e^{-\lambda \underline{F}})$
or 
$\Pi((0,\infty)) = 0$ and $\overline{F}_{\xi^-}>(c/\kappa) (e^{\kappa \underline{F}}-1)$,

-- when $c < 0$ and $\overline{F}<\infty$, if $\Pi((-1, 0)) = 0$ and $\underline F_{\xi}<(c/\lambda) (1-e^{-\lambda \overline{F}})$ 
or 
$\Pi((0,\infty)) = 0$ and $\underline F_{\xi}<(c/\kappa)(e^{\kappa \overline{F}}-1)$.

\end{theorem}

 \section{Ruin probabilities and distributional equations } 
\subsection{Prerequisites}
For the reader convenience we recall the relation the ruin problems with the implicit renewal theory (called also the Kesten--Goldie theory, theory of distributional equations etc.). 
 
\smallskip
First, note that the solution of the non-homogeneous equation (\ref{X}) admits an explicit expression via a stochastic version of the Cauchy formula. In particular, 
$$
X^u_{T^k_n}=e^{V_{T^k_n}}(u-Y_{T^k_n}),
$$ 
where $k\ge 1$, $T^k_n:=T_{kn}$, 
$$
Y_{T^k_n}=-\int_{]0,T^k_n]}e^{-V_{s-}}dP_s=-\sum_{m=1}^n e^{-V_{T^k_{m-1}}}\int_{]T^k_{m-1},T^k_m]}e^{-(V_{s-}-V_{T^k_{m-1}}}dP_s.
$$
Introducing the abbreviations 
$$
M^k_{j}:=e^{-(V_{T^k_{j}}-V_{T^k_{j-1}})}, \qquad Q^k_m=-\int_{]T^k_{m-1},T^k_m]}e^{-(V_{s-}-V_{T^k_{m-1}})}dP_s
$$ 
we rewrite the above formula in the form 
$$
Y_{T^k_n}=\sum_{m=1}^nM^k_1\dots M^k_{m-1}Q^k_m. 
$$
It is easily seen that $\{(M^k_m,Q^k_m)_{m\ge 1}\}$ is an i.i.d. sequence, $\E[|M^k_1|^p=(E[M^1_1|)^p$, and 
$\E[|Q^k_1|^\beta]=\E[|Y_{T_k}|^\beta]<\infty$, cf.
Corollary 2.2 in \cite{EKS}.

The following lemma provides bounds for the ruin probability in terms of the tail of the law of $Y_\infty$: 
\begin{lemma}[{\rm \cite{EKS}, Lemma 3.1}]
\label{G-Paulsen} 
If $Y_{T_n}\to Y_{\infty}$ almost surely as $n\to \infty$, where $Y_{\infty}$ is a finite random variable unbounded from above, 
then for all $u>0$ 
\beq
\label{Paulsen}
\bar G(u)\le\, \Psi(u)
\le {\bar G(u)}/{\bar G(0)},
\eeq 
where $\bar G(u):=\P(Y_\infty>u)$. 
\end{lemma}

Under our assumption the convergence $Y_{T_n}\to Y_{\infty}$ is easy to prove, see \cite{EKS}, Lemma 4.1. Passage to the limit leads to the representation 
\begin{equation}
\label{Y_inf}
 Y_{\infty} = Q^k_1+M^k_1\tilde Y_{\infty},\quad \tilde Y_{\infty} \stackrel{d}{=} Y_{\infty},
\end{equation}
where $\tilde Y_\infty$ is independent of the two-dimensional random variable $(M^k_1,Q^k_1)$. 
This means that we solve the distributional equation 
$Y \stackrel{d}{=}AY+B$, 
 where the law of $(A,B)$ is the law of $(M^k_1,Q^k_1)$.

For this distributional equation we have the following 
\begin{theorem}[{\rm see \cite{EKS}, Theorem 4.2}] 
\label{M}
Suppose that $\beta>0$,
\begin{align}\label{3.3}
\E[M^\beta]=1, \ \ \ \E[M^\beta\,(\ln M)^+]<\infty, \ \ \ \E[|Q|^\beta]<\infty. 
\end{align}
Then $\limsup u^\beta \bar G(u)<\infty$ where $\bar G(u)=\P[Y>u]$. If $Y$ is unbounded from above, then 
$\liminf u^\beta \bar G(u)>0$. 
\end{theorem}
\par
The asymptotic behavior of the ruin probability formulated by (\ref{maina}) follows directly from Lemma \ref{G-Paulsen} and Theorem 
\ref{M} provided that the random variable $Y=Y_\infty$ is unbounded from above. Existing proofs of the latter property are cumbersome and require additional assumptions, see 
\cite{EKS,KPro}. In the present paper we use a rather delicate result, namely, Proposition 2.5.4 from the book \cite{BDM} to get a complement to theorems of the mentioned papers. 

The concept of solution of distributional equation in this book resembles the concept of weak solution in the theory of SDE. In our notations it can be formulated as follows. We are given a measure $m$ on $\bbr^2$. The solution of the distributional equation $Y \stackrel{d}{=}AY+B$ is a probability space $(\Omega,\cF,\P)$ with random variables $(A,B,Y)$ such that $Y$ and $(A,B)$ are independent, the law $\cL(A,B)$ is $m$, and the law $\cL(AY+B)$ coincides with the law $\cL(Y)$. It remains to provide conditions for the existence and uniqueness solution.
 
At first glance, our setting is more specific and resembles the concept of strong solution. Namely, we are given the probability space with an i.i.d. sequence of two-dimensional random variables $(M^k_i,Q^k_i)$ having the same law (which we can denote by $m$). Using moment properties, we construct on this probability space a random variable $\tilde Y_\infty$ independent on $(M^k_1,Q^k_1)$ and having the law which coincides with law of $M^k_1\tilde Y_\infty+Q^k_1$. This gives an idea how to proceed in the abstract setting and get the existence result.

\smallskip 
To formulate Proposition 2.5.4 we need some notations. Associate with the point $(a,b)\in \bbr^2$ the affine transformation of $\mathbb{R}$, namely, the mapping $h(x)= ax+b$. This is a bijection and one can identify $h$ with the point $(a,b)$. 
If $h_i = (a_i b_i)$, $i = 1, 2$, then $h_1 h_2(x) = (a_1 a_2) x + (b_1 + a_1 b_2)$. This property allows us to consider the semigroup of affine transformations of a straight line: $\operatorname{Aff}(\mathbb{R}) = \{h=(a,b) \in \mathbb{R}^2\}$ with the multiplication
$h_1 h_2 := (a_1 a_2, b_1 + a_1 b_2)$.

If $a\neq 1$, then $h$ admits a fixed point $x_0(h)=b/(1-a)$ solving the equation $h(x)=x$. 

Let $\Gamma$ is a subset of $\bbr^2$. Define the semigroup 
$$G(\Gamma): = \{h_1 h_2 \ldots h_n\colon\ h_i \in \Gamma,\ i = 1, 2, \ldots n,\ n \in \mathbb{N}\}.$$ 
 
\begin{proposition}[{\rm \cite{BDM}, Proposition 2.5.4}]
\label{BDM_p254}
Let $Y$ be the solution of distributional equation $Y \stackrel{d}{=}AY + B$. Suppose that there are points $h = (a,b)$ and $h' = (a',b')$ in $G(\operatorname{supp} \cL(A,B))$ such that $0 < a < 1$, $a' > 1$, and 
$$
x_0 (h') := \dfrac{b'}{1 - a'} < x_0 (h): = \dfrac{b}{1-a}.
$$
Then $[x_0(h), \infty) \subset \operatorname{supp} \cL(Y)$.
\end{proposition}

\par
Recall that $\cL(Y_{\infty})$ solves distributional equations with $(A,B)=(M^k_1,Q^k_1)$ and we have a freedom to chose 
the integer parameter $k\ge 1$. The above general proposition implies the following assertion which we shall use in the proof of our results. 

\begin{corollary}\label{coro}
Suppose that there is $k\ge 1$ such that the sets 
$$\{M^k_1>1,\ Q^k_1>0\}, \qquad \{M^k_1<1,\ Q^k_1>0\}$$ are non-null. Then the support of $\cL(Y_\infty)$ is unbounded from above. 
\end{corollary}

In other words, the support of $\cL(Y_\infty)$ is unbounded from above, if there is $k\ge 1$ such that the law of $(M_1^k,Q_1^k)$
charges the sets $\{x>1,\; y>0\}$ and $\{x<1,\; y>0\}$. 

It is worth to notice that  the law of $Y_\infty$ is the solution of various of distributional equation depending on the integer-valued parameter and we are playing with a suitable choice of this parameter and  use  the above proposition in a rather rudimentary form.     

\section{Proofs}

\subsection{Non-life insurance: jumps downwards}

We start from the case where the price process has positive jumps. 
\begin{proposition}
\label{two}
Suppose that $c\! \ge\! 0$, $F_\xi$ charges only $(-\infty,0)$, $\Pi((-1, 0)) = 0$,
$\underline F>0$, and the following condition holds: 

{\bf H.1.} $\overline{F}_{|\xi|}>(c/\lambda) (1-e^{-\lambda \underline F}) $, where $\lambda:=\Pi(h)-a>0$. 

Then the sets $\{M^k_1>1,\ Q^k_1>0\}$ and $ \{M^k_1<1,\ Q^k_1>0\}$ are non-null for any integer $k> 1/(\lambda \underline F)$. 
\end{proposition}

\noindent
{\sl Proof.} 
In the considered case, $\Pi_V((-\infty, 0)) = 0$, $\Pi_V(h)=\Pi(h(\ln (1+x)))$, and $\Pi_V(h) <\infty$. 
The logprice process 
\bean
V&=&
at+ h*(\mu-\nu)+\big (\ln (1+x)-h\big )*\mu=(a-\Pi(h))t+\ln (1+x)*\mu\\
&=&-\lambda t +\ln (1+x)*\mu
=
-\lambda t + x*\mu^V, 
\eean
where $\lambda :=\Pi (h)-a> 0$, otherwise Assumption \ref{ass1} is violated. 

\smallskip

Fix $T<\infty$ and $\e\in (0,1]$. Let us show that $V$ admits the representation 
$$
V=-\lambda t+ Z^1+Z^2+Z^3,
$$
where the increasing processes $Z^i$ are independent, $Z^2$ and $Z^3 $ are compound Poisson, 
the set $\{Z^1\le \e\}$ is non-null and the range of jumps $[\hat p,p]$ of $Z^2$ is arbitrarily narrow. 
To this aim we consider two cases. 

${\bf a.}$ The point $\underline \Pi_V=:p$ is an atom of $\Pi_V$. To meet the above requirements we take $Z^1:=0$, $Z^2:=xI_{\{p\}}*\mu$, and $Z^3:=xI_{(p,\infty)}*\mu$. 
\smallskip

${\bf b.}$ The point $\underline \Pi_V$ is not an atom of $\Pi_V$. Define the process $Z^1:=xI_{(\underline \Pi_V,\bar p)}*\mu$, where $\bar p\in (\underline{\Pi}_V,\overline{\Pi}_V)$, sufficiently close to $\underline \Pi_V$, is such that $\P(Z^1_T\le \e)>0$. 
Such a choice is possible since $\E\big [ I_{[\underline{\Pi}_V, \hat p)}h*\mu^V_T\big]=\Pi_V \big(I_{[\underline{\Pi}_V,\hat p)}h\big)T\downarrow 0$ as $\hat p \downarrow \underline{\Pi}_V $
 and, in virtue of the Chebyshev inequality, $\P[Z^1_T > \e] < 1$ when $\hat p$ is close to $\underline \Pi_V$. 
We can take $Z^2:=xI_{[\hat p, p]}*\mu$ and $Z^3:=xI_{(p,\infty)}*\mu$ with $p>\hat p$ arbitrary close to $\hat p$.

\bigskip
We associate with the compound Poisson process $Z^2$ the counting Poisson process $N^2:=I_{[\hat p,p]}*\mu^V_t$ and define, for arbitrary 
$T',T$ such that $0\le T'<T<\infty$, the sets 
$C_n(0,T'):=\{ N^2_{T'}- N^2_{0}=n\}$ and $C_0(T',T):=\{ N^2_{T}- N^2_{T'}=0\}$. In other words, we consider the sets where $Z^2$ has exactly $n$ jumps on the interval $(0,T']$ and has no jumps on $(T',T]$. 
They are non-null as well as the set 
$$
A(T,n) := \{ Z^3_T = 0 \} \cap \left\{ Z^1_T \le \e \right \} \cap C_n(0,T')\cap C_0(T',T).
$$ 
Take $\Delta>0$ and define the set $B_k := \bigcap_{1 \le j \le k} \{T_j - T_{j - 1} \in [\underline F, \underline F + \Delta ) \}$. It is also non-null and on it the arrival times $T_j\in [j\underline F, j\underline F + j\Delta )$. 

Put $T=T(k):= k\underline F + k\Delta$. On the non-null set $A(T(k),n)\cap B_k$ the process 
$V$ coincides with $Z^3$ on the time interval $[0,T(k)]$ while $Z^1$ remains below $\e$ at this interval. Moreover, for any $T'\in (0,\underline F)$
\bea
\nonumber
\int_{[0,{T_k}]} e^{-V_{s}} ds
&\le& 
\int_{[0,T']} e^{\lambda s} ds +e^{-n\hat p}\int_{(T',T]} e^{\lambda s} ds \\
\label{int-l}
&=&\frac {1}
{\lambda}
 (e^{\lambda T'}-1)+\frac {1}
{\lambda}e^{-n \hat p} (e^{\lambda T}-e^{\lambda T'}) 
\eea
and, for every $j\le k$, 
\bea
\label{e-j-g} 
e^{-V_{T_{j}}} &=& e^{ \lambda T_j-Z^1_{T_j}-Z^2_{T_j}} \ge e^{ \lambda j \underline F-\e -Z^2_{T_j}}\ge e^{\lambda j \underline F-\e -np }, \\
\label{e-j-l}
e^{-V_{T_{j}}} &=& e^{ \lambda {T_j}-Z^1_{T_j}-Z^2_{T_j}} \le e^{\lambda j (\underline F+\Delta)}e^{-Z^2_{T_j}}\le e^{\lambda j (\underline F+\Delta)-n \hat p }. 
\eea
In particular, for $M^k_1:=e^{-V_{T_{k}}}$ there are the bounds
\beq
\label{M-k}
M^k_1\ge e^{\lambda k \underline F-\e -np }, \qquad
M^k_1\le e^{\lambda k (\underline F+\Delta)-n \hat p }. 
\eeq

\medskip
Suppose $\overline{F}_{|\xi|}<\infty$ (the opposite case is covered by Theorem \ref{old1}). 
Assume first that $\underline F_{|\xi|}< \overline{F}_{|\xi|}$.
Take $\delta \in (0,\overline{F}_{|\xi|}-\underline F_{|\xi|})$.
Then $ \Gamma_k(\delta):=\cap_{j=1}^k\{|\xi_j|> \overline{F}_{|\xi |}- \delta\} $ is a non-null set as well 
as the set $A(T(k),n)\cap B_k \cap \Gamma_k(\delta)$. On the latter we have the bound 
$$
Q_1^k> -\frac c \lambda \big[
 (e^{\lambda T'}-1)+e^{-n \hat p} (e^{\lambda T}-e^{\lambda T'})\big]+ ( \overline{F}_{|\xi|}-\delta)\frac {e^{\lambda \underline F} }{e^{\lambda \underline F}-1} e^{- \e -n p} (e^{\lambda k \underline F}- 1), 
$$
where $T'\in (0,\underline F)$. 
Regrouping terms in the right-hand side, we rewrite the above inequality in the form 
$$
Q_1^k> \left(\gamma +R_1+R_2\right)e^{ -n p} (e^{\lambda k\underline F}- 1),
$$
where 
$$
\gamma:=\overline{F}_{|\xi|}\frac {e^{\lambda \underline F} }{e^{\lambda \underline F}-1}-\frac c\lambda>0
$$
by virtue of {\bf H.1.} The term 
$$
R_1:=(e^{-\e }-1)\overline{F}_{|\xi|}\frac {e^{\lambda \underline F} }{e^{\lambda \underline F}-1}- \delta e^{-\e }\frac {e^{\lambda \underline F} }{e^{\lambda \underline F}-1} - \frac c\lambda 
\frac {e^{\lambda T'}-1}{e^{ \lambda k \underline F}- 1}e^{np}>-\gamma/2
$$
for any $n$ by a choice of sufficiently small $\e$, $\delta$ and $T'=T'(n,p)$, and the term 
\bean
R_2:= \frac c\lambda - \frac c\lambda e^{n (p-\hat p)} \frac{e^{\lambda T(k)}-e^{\lambda T'}}{e^{ \lambda k \underline F}- 1}
&= & \frac c\lambda \left[1
- e^{n (p-\hat p)} \frac{e^{\lambda k (\underline F+\Delta)}-1}{e^{ \lambda k \underline F}- 1}\right ]\\
&& + \frac c\lambda e^{n (p-\hat p)} \frac{e^{\lambda T'}-1}{e^{ \lambda k \underline F}- 1}. 
\eean

Take $n=0$ and sufficiently small $ \Delta$ and $T'$ to ensure that 
$R_2> -\gamma/2$. 
 Then the set $\{M^k_1>1,\ Q^k_1>0\}$ contains the non-null set $A(T(k),0)\cap B_k \cap \Gamma_k(\delta)$. 
 
 The set $\{M^k_1< 1,\ Q^k_1>0\}$ will contain 
 the 
 set $A(T(k),n)\cap B_k \cap \Gamma_k(\delta)$, if we take an arbitrary integer $n\ge e^{\lambda k (\underline F+1)}/\hat p$, and sufficiently small values of $\Delta$ and $p-\hat p$ to ensure that 
$$
e^{n (p-\hat p)} \frac{e^{\lambda k (\underline F+\Delta)}-1}{e^{ \lambda k \underline F}- 1}<1+\frac {\gamma c}{2\lambda}. 
$$ 
 
 The case $\underline F_{|\xi|}= \overline{F}_{|\xi|}$, i.e. $\xi$ is a constant, is easy. On the set 
 $A(T(k),n)\cap B_k$ we have the same bound for $Q_1^k$ but now with $\delta=0$ and we conclude as above. \fdem

\medskip

\begin{proposition}
\label{three}
Suppose that $c\! \ge\! 0$, $F_\xi$ charges only $(-\infty,0)$, $\Pi((0,\infty)) = 0$, 
$\underline F>0$, and the following condition holds: 

{\bf H.2.} $\overline{F}_{|\xi|}>(c/\kappa) (e^{\kappa \underline F}-1) $ where $\kappa:=a-\Pi(h)>0$. 

Then the sets $\{M^k_1>1,\ Q^k_1>0\}$ and $ \{M^k_1<1,\ Q^k_1>0\}$ are non-null for any integer $k> 1/(\kappa \underline F)$. 
\end{proposition}
{\sl Proof.}
Take $T<\infty$ and $\e>0$. 
By analogy with the arguments of previous subsection we split the range of jumps of $V$ by appropriate points 
(negative now) 
and represent $V$ as the sum of independent processes
\beq
\label{pi-neg}
V=\kappa t+xI_{(-\infty,-p)}*\mu^V+xI_{[-p,-\hat p]}*\mu^V+xI_{(-\hat p,0)}*\mu^V= \kappa t+Z^3+Z^2+Z^1,
\eeq
where 
$Z^3:=xI_{(-\infty,-p)}*\mu^V$ and $Z^2:=I_{[-p,-\hat p]}*\mu^V$ are compound Poisson processes, 
$Z^1:=xI_{[-\hat p,0)}*\mu^V$, $p\ge \hat p>0$ and $Z^1_T>-\e$ on the non-null set. If $\overline{\Pi}_V$ is the atom of $\Pi_V$, we take $p=\hat p=-\overline{\Pi}_V$ and $Z^1\equiv 0$. The difference $p-\hat p$ can be chosen arbitrarily small. 

Retaining the notations of the previous subsection for $C_n$ and $C_0$, we obtain, taking $T'\in (0,\underline F)$, that on the non-null set 
$$
A(T,n)\cap B_k := \{ Z^3_T = 0 \} \cap \left\{ Z^1_T >-\e \right \} \cap C_n(0,T')\cap C_0(T',T)\cap B_k,
$$
the following inequalities hold: 
\bean
\label{int-l-}
\int_{[0,{T_k}]} e^{-V_{s}} ds
&\le& 
e^{np+\e }T'
+e^{n p}\int_{[T',T(k)]} e^{-\kappa s} ds \\
&&
=
e^{np+\e }T'
 +\frac {1}
{\kappa}e^{n p+\e } (e^{-\kappa T'}-e^{-\kappa T(k)}),\\ 
e^{-V_{T_{j}}} 
 &\ge& 
e^{-\kappa j (\underline F+\Delta)+n\hat p }, \qquad
e^{-V_{T_{j}}} \le e^{-\kappa j \underline F+\e +n p }, \qquad j\le k. 
\eean
In particular, for $M^k_1:=e^{-V_{T_{k}}}$ there are the bounds
$$
M^k_1\ge e^{-\kappa k (\underline F+\Delta)+n\hat p }, \qquad
M^k_1\le e^{-\kappa k \underline F+\e +n p }. 
$$

\smallskip
Suppose that $\overline{F}_{|\xi|}<\infty$. 
Assume that $\underline F_{|\xi|}< \overline{F}_{|\xi|}$.
Take $\delta \in (0,\overline{F}_{|\xi|}-\underline F_{|\xi|})$ and put
 $$
 \Gamma_k(\delta):=\cap_{j=1}^k\{|\xi_j|> \overline{F}_{|\xi |}- \delta\}.
 $$
 On the non-null set 
$A(T(k),n)\cap B_k \cap \Gamma_k(\delta)$ we have the bound  
\bean
Q_1^k&>& -c e^{np+\e }T'
 -\frac {c}
{\kappa}e^{n p+\e } (e^{-\kappa T'}-e^{- \kappa k (\underline F+\Delta)})\\
&&
 + (\overline{F}_{|\xi |}- \delta) e^{n\hat p } \frac {e^{-\kappa (\underline F+\Delta)} }{1-e^{-\kappa (\underline F+\Delta)}} (1-e^{- \kappa k (\underline F+\Delta)}). 
\eean
 Regrouping terms in the right-hand side we rewrite the above inequality 
in the form 
$$
Q_1^k> \left( \gamma +R_1+R_2\right)e^{ n \hat p} (1-e^{- \kappa k( \underline F+\Delta)}),
$$
where 
$$
\gamma:=\overline{F}_{|\xi|}\frac {e^{-\kappa \underline F} }{1-e^{-\kappa \underline F}}-\frac c\kappa>0
$$
by virtue of {\bf H.2.} and 
\bean
R_1&:=&\overline{F}_{|\xi|}\left[\frac {e^{-\kappa (\underline F+\Delta)} }{1-e^{-\kappa (\underline F+\Delta)}}-\frac {e^{\kappa \underline F} }{e^{\kappa \underline F}-1}\right]- \delta \frac {e^{-\kappa (\underline F+\Delta)} }{1-e^{-\kappa (\underline F+\Delta)}}\\
&& - 
\frac {ce^{n(p-\hat p)}}{1-e^{- \kappa k( \underline F+\Delta)}}\left(e^{\e } T'+\frac 1\kappa (e^{\kappa T'}-1)\right),\\
R_2&:=& \frac c \kappa(1-e^{n (p-p')}). 
\eean

Take $n=0$ and sufficiently small $\delta, \Delta, T'$ to ensure that 
$R_1> -\gamma$. 
 Then the set $\{M^k_1<1,\ Q^k_1>0\}$ contains the non-null set $A(T(k),0)\cap B_k \cap \Gamma_k(\delta)$. 
 
 Finally, take an integer $n\ge e^{\kappa k (\underline F+1)}/\hat p$. Chose $\delta$, $\Delta$, $T'$, and 
 $p-\hat p$
 small enough to ensure that $R_1>-\tilde \gamma/2$ and $R_2>-\tilde \gamma/2$. 
 Then $\{M^k_1> 1,\ Q^k_1>0\}$
 contains the non-null set $A(T(k),n)\cap B_k \cap \Gamma_k(\delta)$. \fdem 

\begin{remark}
In the formulations of the above propositions we can replace $|\xi|$ by $\xi^-$ (recall that $\overline{F}_\xi\le 0$). 
In such a modified form the conclusions hold also for the mixed model, that is, in the case of two-sided jumps of the business process, where $\overline{F}_{\xi^+}>0$ as well as $\overline{F}_{\xi^-}>0$. The proofs need only minor changes. 
\end{remark} 

\subsection{Annuity payments: jumps upwards}
Again we consider first the case where the price process has positive jumps. 
\begin{proposition}
\label{four}
Suppose that $c\! < \! 0$, $\underline F_\xi\ge 0$, $\Pi((-1, 0)) =0$, 
$\overline{F}<\infty$, and the following condition holds:

{\bf H.3.} $\underline F_{\xi}<(c/\lambda) (1-e^{-\lambda \overline{F}}) $, where $\lambda:=\Pi(h)-a>0$. 

Then the sets $\{M^k_1>1,\ Q^k_1>0\}$ and $ \{M^k_1<1,\ Q^k_1>0\}$ are non-null for any integer $k> 2/(\lambda \overline{F})$. 
\end{proposition}
{\sl Proof.} 
We proceed as in the proof of Proposition \ref{two} representing the logprice process in the form $
V=-\lambda t+ Z^1+Z^2+Z^3$, where the increasing jump processes $Z^i\ge 0$. Since $c$ is now negative, we 
consider sets where the integral with respect to $ds$ has larger values comparatively to the sum of jumps. 
To this aim, we introduce the non-null sets 
\bean
A_k(\Delta,n)& := &\{ Z^3_T = 0 \} \cap \left\{ Z^1_{T} \le \e \right \} \cap C_0(0,(k-1)\overline{F})\cap C_n(k\overline{F}_\Delta,k\overline{F}), \\
 B_k(\Delta)&:=&\cap_{j=1}^k\{T_j-T_{j-1}\in (\overline{F}_\Delta,\overline{F}]\},
\eean
where $T=k\overline{F}$, $\Delta\in (0,\overline{F}/2)$, $\overline{F}_\Delta:=\overline{F}-\Delta$, and $\e\in (0,1]$. On the intersection of these sets, $T_j\in (j\overline{F} _\Delta,j \overline{F}]$ for all $j\le k$ and
\bean
\label{int-g}
\int_{[0,{T_k}]} e^{-\lambda V_{s}} ds
\ge 
\int_{[0,k\overline{F}_\Delta]} e^{\lambda s} ds 
=
\frac 1\lambda 
(e^{\lambda k\overline{F}_\Delta})-1).
\eean
Moreover, $Z^2_{T_j}=0$ for $j<k$ and, therefore,
\bean
\label{bounds-j} 
e^{-V_{T_{j}}} \ge 
 e^{\lambda j \overline{F}_\Delta-\e}, \qquad
e^{-V_{T_{j}}} \le 
 e^{\lambda j \overline{F}}. 
\eean
Since the process $Z^2$ has $n$ jumps on the interval $(k\overline{F}_\Delta,k\overline{F}
]$, the bounds are different for $j=k$: 
\bean
M^k_1\ge e^{\lambda k \overline{F}_\Delta-\e-np }, \qquad
M^k_1\le e^{\lambda k \overline{F}-n \hat p }. 
\eean

Let $\Gamma(\delta):=\cap_{j=1}^k\{\xi_j\le \underline F_\xi+\delta\}$. Then on the set $A_k(\Delta,n)\cap B_k(\Delta)\cap \Gamma(\delta) $ we have the bound 
$$
\sum_{j=1}^ke^{-V_{T_j}}\xi_j\le (\underline F_\xi + \delta) \left(\frac {e^{\lambda \overline{F}} }{e^{\lambda \overline{F}}-1} (e^{\lambda k\overline{F}}-1)+e^{\lambda k \overline{F}-n \hat p }-e^{\lambda k \overline{F}}\right). 
$$
Combining it with (\ref{int-g}) we get that 
$$
Q_1^k\ge \frac{|c|}\lambda 
(e^{\lambda k \overline{F}_\Delta}-1)-(\underline F_\xi + \delta) \left(\frac {e^{\lambda \overline{F}} }{e^{\lambda \overline{F}}-1} (e^{\lambda k\overline{F}}-1)+e^{\lambda k \overline{F}}(e^{-n \hat p }-1)\right). 
$$
Regrouping terms in the right-hand we rewrite this inequality in the form 
$$
Q_1^k\ge \big(\gamma+R_1+R_2\big) \big (e^{\lambda k\overline{F}}-1\big), 
$$
where the leading term 
$$
\gamma:=\frac{|c|}{\lambda}-\underline F_\xi \frac {e^{\lambda \overline{F}}}{e^{\lambda \overline{F}}-1}>0
$$
due to the assumption {\bf H.3.} 
 The term
$$
R_1:=\frac{|c|}\lambda \left(\frac 
{e^{\lambda k \overline{F}_\Delta} -1}{e^{\lambda k\overline{F}} -1} - 1\right) - \delta \left(\frac {e^{\lambda \overline{F}} }{e^{\lambda \overline{F}}-1} +\frac{e^{\lambda k \overline{F}}}{e^{\lambda k\overline{F}} -1} (e^{-n \hat p }-1)\right)>-\gamma/2 
$$
independently on $n$, when $\Delta$ and $\delta$ are suffitiently small. The term 
$$
R_2:=\underline F_\xi \frac{e^{\lambda k \overline{F}}}{e^{\lambda k\overline{F}} -1}(1-e^{-n\hat p})> 0. 
$$ 

With our choice of $k$, $\Delta$, and $\delta$ the set $\{M^k_1> 1,\ Q^k_1>0\}$ contains the non-null set $A_k(\Delta,0)\cap B_k(\Delta) \cap \Gamma_k(\delta)$. 
The set $\{M^k_1<1,\ Q^k_1>0\}$ is also non-null since it contains, for sufficiently small $\Delta$ and $\delta$, the set $A_k(\Delta,n)\cap B_k(\Delta) \cap \Gamma_k(\delta)$, when the integer $n>\lambda k\overline{F}/\hat p$. \fdem 
\begin{proposition}
\label{five}
Suppose that $c\! < \! 0$, $\underline F_\xi\ge 0$, $\Pi((0,\infty)) =0$, 
$\overline{F}<\infty$, and the following condition holds:

{\bf H.4.} $\underline F_{\xi}<(c/\kappa)(e^{\kappa \overline{F}}-1) $, where $\kappa:=a-\Pi(h)>0$. 

Then the sets $\{M^k_1>1,\ Q^k_1>0\}$ and $ \{M^k_1<1,\ Q^k_1>0\}$ are non-null for all sufficiently large integers $k$. 
\end{proposition}
{\sl Proof.} In the considered case the logprice $V$ decreases by jumps and increases continuously with 
the drift coefficient $a-\Pi(h)>0$. As in Proposition \ref{three}, we use the representation 
$V= \kappa t+Z^3+Z^2+Z^1$
with negative decreasing processes $Z^i$, $i=1,2,3$, where the central role plays the compound Poisson process $Z^2$ with jumps in an appropriately chosen range. 
As in Proposition \ref{four} we are interested in sets where the integral with respect to $ds$ (the annuity payments) has larger values with respect to the sum of gains $\xi_i$. We introduce the non-null sets
\bean
A_k(\Delta,\e,n)& := &\{ Z^3_T = 0 \} \cap \left\{ Z^1_{T} \ge -\e \right \} \cap C_0(0,(k-1)\overline{F})\cap C_n(k\overline{F}_\Delta,k\overline{F}), \\
 B_k(\Delta)&:=&\cap_{j=1}^k\{T_j-T_{j-1}\in (\overline{F}_\Delta,\overline{F}]\},
\eean
where $\Delta\in (0,\overline{F}/2)$ is a small number and $\overline{F}_\Delta:=\overline{F}-\Delta$. On the intersection of these sets, the renewal times $T_j\in (j\overline{F} _\Delta,j \overline{F}]$ for all $j\le k$ and 
\beq
\label{int-g1}
\int_{[0,{T}]} e^{-\kappa V_{s}} ds
\ge 
\int_{[0,k\overline{F}_\Delta]} e^{-\kappa s} ds 
=
\frac 1\kappa 
(1-e^{- \kappa k \overline{F}_\Delta}). 
\eeq
Moreover, $Z^2_{T_j}=0$ for $j<k$ and, therefore, 
\bean
e^{-V_{T_{j}}} \ge
e^{-\kappa j \overline{F}}, \qquad
e^{-V_{T_{j}}} \le 
 e^{-\kappa j \overline{F}_\Delta+\e}. 
\eean
Since the process $Z^2$ has $n$ jumps on the interval $(k\overline{F}_\Delta,k\overline{F})$, for $j=k$ we have the bounds
\bean
M^k_1\ge e^{- \kappa k\overline{F}+n\hat p }, \qquad
M^k_1\le e^{- \kappa k\overline{F}_\Delta +\e + n p }. 
\eean

Let $\Gamma(\delta):=\cap_{j=1}^k\{\xi_j\le \underline F_\xi+\delta\}$. Then on the set $A_k(\Delta,\e,n)\cap B_k(\Delta)\cap \Gamma(\delta) $
$$
\sum_{j=1}^k e^{-V_{T_j}}\xi_j\le (\underline F_\xi + \delta)e^\e \left(\frac {e^{-\kappa \overline{F}_\Delta} }{1-e^{-\kappa 
\overline{F}_\Delta}} (1-e^{-\kappa k \overline{F}_\Delta})+e^{-\kappa k \overline{F}_\Delta}(e^{n p }-1)\right).
$$
Combining this bound with (\ref{int-g1}) we get that 
$$
Q_1^k\ge \frac{|c|}\kappa 
(1-e^{-\kappa k\overline{F}_\Delta})-(\underline F_\xi + \delta)e^\e\Bigg(\frac {e^{-\kappa \overline{F}_\Delta} }{1-e^{-\kappa \overline{F}_\Delta}} (1-e^{- \kappa k \overline{F}_\Delta})+e^{- \kappa k \overline{F}_\Delta}(e^{n p }-1)\Bigg).
$$
Let us denote the expression in parentheses by $g(\Delta, np)$. Regrouping terms in the right-hand we rewrite the above inequality in the form 
$$
Q_1^k\ge \big(\gamma+R_1+R_2\big)e^{n\hat p+\e} \big (1-e^{-k \kappa \overline{F}}\big), 
$$
where the leading term 
$$
\gamma:=\frac{|c|}{\kappa}-\underline F_\xi \frac {e^{-\kappa \overline{F}}}{1-e^{-\kappa \overline{F}}}>0
$$
due to the assumption ${\bf H.3.}$ The term 
$$
R_1=R_1(\Delta, \delta, \e,n):=\frac{|c|}\kappa e^{-n\hat p-\e} \left ( \frac 
{1-e^{- \kappa k \overline{F}_\Delta} }{1-e^{-\kappa k \overline{F}}} - 1\right)- \delta e^{-n\hat p}\frac{g(\Delta, np)}{1-e^{- \kappa k \overline{F}}}
$$
dominates $-\gamma/2$ for sufficiently small $\Delta$ and $\e$. The term 
$$
R_2=R_2(\Delta, \e, n):= \frac{\underline F_\xi}{1-e^{- \kappa k\overline{F}}} e^{-n\hat p} \left(e^{-\e}\frac {e^{-\kappa \overline{F}}}{1-e^{-\kappa \overline{F}}}-g(\Delta,np)\right ). 
$$

Note that the function $g$ is continuous and 
$$
g(\Delta, 0)\to \frac {e^{-\kappa \overline{F}}}{1-e^{-\kappa \overline{F}}}(1-e^{- \kappa k\overline{F}}) \qquad \hbox{as }\ \Delta\to 0. 
$$
It follows that if we take $n=0$, we can find sufficiently small $\Delta,\delta,\e$ to ensure that 
$R_1(\Delta,\delta,\e,0)>-\gamma$ and $R_2(\Delta,0)>0$. 
Therefore, for any integer $k> 2/(\kappa \overline{F})$ the set $\{M^k_1< 1,\ Q^k_1>0\}$ contains the non-null set 
$A_k(\Delta,\e,0)\cap B_k(\Delta) \cap \hat \Gamma_k(\delta)$. 
 
To complete the proof, we enlarge, in the need, the integer $k$ to ensure the bound 
$$
\frac{\underline F_\xi e^{-\kappa k\overline{F}}}{1-e^{- \kappa k\overline{F}}}<\gamma/8. 
$$ 
Take arbitrary $\hat p>0$, chose $n$ such that $n\hat p>\kappa k \overline{F}$ and $p>\hat p$ sufficiently close to $\hat p$ to satisfy the inequality $e^{n(p-\hat p)}\le 2$. Note that 
$$
g(\Delta, np)\to \frac {e^{-\kappa \overline{F}}}{1-e^{-\kappa \overline{F}}}(1-e^{- \kappa k\overline{F}})+e^{-\kappa k\overline{F}} (e^{np}-1)\qquad \hbox{as }\ \Delta\to 0, 
$$
and 
$$
R_2(\Delta, \e, n)\to \frac{\underline F_\xi e^{- \kappa k\overline{F}}}{1-e^{- \kappa k\overline{F}}} 
\left(\frac {e^{-n\hat p-\e} }{1-e^{-k \overline{F}}}-e^{n(p-\hat p)}\right)>-\gamma/4.
$$ 
Thus, for such a choice of parameters, $R_2(\Delta, \e, n)>-\gamma/2$ when $\Delta$ is sufficiently small. It follows that the set $\{M^k_1>1,\ Q^k_1>0\}$ is non-null, since it contains for a suitably chosen $k$, $n$, $\hat p$, $p$ and sufficiently small values $\Delta$ and $\delta$ the non-null set $A_k(\Delta,n)\cap B_k(\Delta) \cap \Gamma_k(\delta)$. \fdem 

\section{Comments on the mixed model}

Let us consider the case where both $\overline{F}_{{\color{red}\xi^+}}$ and $\overline{F}_{{\color{red}\xi^-}}$ are strictly positive. Suppose that $c\ge 0$. An inspection of the proofs of Propositions \ref{two} and \ref{three} shows that arguing on the sets where $\xi_j<0$ we get the following 
\begin{proposition}
\label{two-three}
Suppose that $c\! \ge\! 0$ and $\underline F>0$. Then the sets $\{M^k_1>1,\ Q^k_1>0\}$ and $ \{M^k_1<1,\ Q^k_1>0\}$ are non-null for any $k> 1/(\lambda \underline F)$ in the following cases:

$(a)$ $\Pi((-1, 0)) = 0$ and $\overline{F}_{\xi^-}>(c/\lambda) (1-e^{-\lambda \underline F}) $ where $\lambda:=\Pi(h)-a>0$; 
 
$(b)$ $\Pi((0,\infty)) = 0$ and $\overline{F}_{\xi^-}>(c/\kappa) (e^{\kappa \underline F}-1) $ where $\kappa:=a-\Pi(h)>0$. 
\end{proposition}

An inspection of the proofs of Propositions \ref{four} and \ref{five} leads to 

\begin{proposition}
\label{four-five}
Suppose that $c< 0$ and $\overline{F}<0$. Then the sets $\{M^k_1>1,\ Q^k_1>0\}$ and $ \{M^k_1<1,\ Q^k_1>0\}$ are non-null for any sufficiently large $k$ in the following cases:

$(c)$ $\Pi((-1, 0)) = 0$ and $\underline F_{\xi^+}<(c/\lambda) (1-e^{-\lambda \overline{F}}) $, where $\lambda:=\Pi(h)-a>0$; 
 
$(d)$ $\Pi((0,\infty)) = 0$ and $\underline F_{\xi^+}<(c/\kappa)(e^{\kappa \overline{F}}-1) $, where $\kappa:=a - \Pi(h)>0$.
\end{proposition}

\smallskip
{\bf Acknowledgements.}

Platon Promyslov expresses his sincere gratitude to the Mathematical laboratory of Besan\c con of the University of Franche-Comté for the opportunity to carry out his scientific research as a visitor under the Ostrogradsky scholarship program.

\section*{Competing interests}
The authors declare no competing interests.

\end{document}